\documentclass[preprint,review,10pt]{elsarticle}
\usepackage{amsfonts}
\usepackage{amsmath,amssymb}
\usepackage{graphicx}
\usepackage{setspace}
\usepackage[left=2cm,top=2cm,right=2cm]{geometry}

\newtheorem{lemma}{Lemma}

\newtheorem{theorem}{Theorem}
\numberwithin{equation}{section}
\journal{FILOMAT}
\begin{document}
\title{Some Approximation Properties of New Families of Positive Linear Operators}
\author{Prashantkumar Patel}
\ead{prashant225@gmail.com}
\address{Department of Mathematics,
St. Xavier's College(Autonomous), Ahmedabad-380 009 (Gujarat),
India}
%\fntext{Corresponding authors}
\begin{abstract}
In the present article, we propose the new class positive linear
operators, which discrete type depending on a real parameters. These operators are similar to Jain operators but its approximation properties are different then Jain operators. Theorems of degree of approximation, direct results, Voronovskaya Asymptotic formula and statistical convergence are discussed.
\end{abstract}

\maketitle
\section{Introduction}
During the last decade two types of generalizations of the
classical Poisson, binomial and negative binomial distribution,
useful in biology, ecology and medicine  have been introduced by
considering the two basic forms of Lagrange series
\begin{equation}\label{18.1.eq1}
\phi(z) =\phi(0)
+\sum_{k=1}^{\infty}\frac{1}{k!}\left[\frac{d^{k-1}}{dz^{k-1}}(f(z))^k\phi'(z)\right]_{z=0}
\left(\frac{z}{f(z)}\right)^k;
\end{equation}
\begin{equation}\label{18.1.eq2}
\phi(z)\left[1-\frac{z}{f(z)}\frac{df(z)}{dz}\right]^{-1} =
\sum_{k=0}^{\infty}\frac{1}{k!}\left[\frac{d^{k}}{dz^{k}}(f(z))^k\phi(z)\right]_{z=0}
\left(\frac{z}{f(z)}\right)^k \end{equation} and expanding
suitable function $\phi(z)$ into powers of $\displaystyle
\frac{z}{f(z)}$ for suitably chosen function $f(z)$.\\
By putting $\phi(z)=e^{\alpha z}$ and $f(z)=e^{\beta z}$ in
formulas (\ref{18.1.eq1}) and (\ref{18.1.eq2}), we achieved
Generalized Poisson Distribution (GPD) studied by Consul and Jain
\cite{consul1973generalization,consul1973some} and Linear Function
Poisson Distribution(LFPD) was introduced and studied by Jain
\cite{jain1975linear}. \\
\indent Since 1912, Bernstein Polynomial and its various
generalization have been studied by Bernstein
\cite{Berstien1912Demo}, Sz\'{a}sz \cite{szasz1950generalization},
Meyer-Konig and Zeller \cite{meyer1960bernsteinsche}, Cheney and
Sharma \cite{cheney1964generalization}, Stancu
\cite{stancu1968approximation}. Bernstein polynomials are based on
binomial and negative binomial distributions. In 1941, Sz\'{a}sz
and Mirakyan \cite{mirakyan1941approximation} have introduced
operator using the Poisson distribution.  %After that various
%generalizations and approximation properties was discussed by many
%researcher \cite{}. 
We mention that rate of convergence developed
by Rempulska and Walczak \cite{rempulska2001approximation},
asymptotic expansion introduced by Abel \textit{et al.}
\cite{abel2007asymptotic}. In 1976, May \cite{may1976saturation}
showed that the Baskakov operators can
reduces to the Sz\'{a}sz-Mirakyan operators.\\
\indent The Lagrange formula \eqref{18.1.eq1} was used to
established the Jain operators \cite{jain1972approximation}, which
are as follows
\begin{equation}\label{18.1.eq3}
 J_{n}^{[\beta]}(f,x) =  \sum_{k=0}^{\infty} \omega_{\beta}(k,n x)
 f\left(\frac{k}{n}\right),
 \end{equation}
 where $\displaystyle \omega_{\beta}(k,n x) = n x( n x  + k\beta)^{k-1} \frac{e^{-(n
 x+k\beta)}}{k!}$; $0\leq \beta< 1$ and $f\in C[0,\infty)$.\\
\indent The operators (\ref{18.1.eq3}) are generalization of well
known the Sz\`{a}sz-Mirakyan operators. In 2013, Agratini
\cite{agratini2013approximation} discussed the relation between
the local smoothness of function and local approximation. Also,
the degree of approximation and the statistical convergence of the
sequence \eqref{18.1.eq3} was studied in
\cite{agratini2013approximation}. We mention that a
Kantorovich-type extension of the Jain operators was given in
\cite{umar1985approximation}. Additionally, the Durrmeyer type
generalization of the Jain operators was established in
\cite{tarabie2012jain,mishrasome2013,Patelmishra20131,agratini2014approximation}.
The Jain operators was also developed in two variables in
\cite{farcasapproximation}. The Jain type variant of Lupas operators \cite{lupas1995approximation} was studied by Patel and Mishra in \cite{patel2015onnew}.  Due to their properties, the operators
$J_{n}^{[\beta]}$ and $J_{n}^{[0]}$ have been intensively studied
by many mathematicians. Thus, in our opinion, the class $P_{n}^{[\beta]}$ defined in \eqref{18.2.eq1} should deeper investigate.\\
\indent In this manuscript, we use Lagrange formula
\eqref{18.1.eq2} to establish new sequence of positive linear
operators. The approximation properties establish in this
manuscript are different then the Jain operators \eqref{18.1.eq3}.
Local approximation properties, the rate of convergence, weighted
approximation, asymptotic formula and  statistical convergence are
investigated for the sequence of the operators \eqref{18.2.eq1}.\\
For $0<\alpha<\infty$ and $|\beta| <1$, proceed by setting $
\phi(z) = e^{\alpha z}$~  and  ~$~f(z) = e^{\beta z},$ in Lagrange
formula \eqref{18.1.eq2}, we get
\begin{eqnarray*}
e^{\alpha z}\left[1-\beta z\right]^{-1} &=&
\sum_{k=0}^{\infty}\frac{1}{k!}(\alpha+k\beta)^k
\left(\frac{z}{e^{\beta z}}\right)^k.
\end{eqnarray*}
Therefore, we shall have
\begin{equation} e^{\alpha z}  = (1-\beta z)\sum_{k=0 }^{\infty}  (\alpha + \beta k )^{k}  \frac{u^k}{k!};~~~~~u=ze^{-\beta z}, \end{equation}
where $z$ and $u$ are sufficiently small such that $|\beta u | <
a^{-1}$ and $| \beta z|<1$.\\
By taking z = 1, we have
\begin{equation}\label{18.1.eq4}
1  = (1-\beta)\sum_{k=0 }^{\infty} \frac{1}{k!} (\alpha + \beta k
)^{k}e^{-(\alpha + \beta k)}.
\end{equation}
Define
\begin{equation}\label{18.1.eq5}\mathbf{p}_{\beta}(k,\alpha)= (1-\beta)\frac{1}{k!} (\alpha + \beta k
)^{k}e^{-(\alpha + \beta k)}.\end{equation} Now, form equality
\eqref{18.1.eq4}, we can write

$$\sum_{k=0}^{\infty}\mathbf{p}_{\beta}(k,\alpha)=1,$$
for $0<\alpha<\infty$ and $|\beta|<1$.
\section{Construction of the Operators}
We may now define the operator as
\begin{equation}\label{18.2.eq1}
P_{n}^{[\beta]}(f,x) = \sum_{k=0}^{\infty} \mathbf{p}_{\beta}(k,n
x) f\left(\frac{k}{n}\right),
\end{equation}
where $0\leq \beta < 1$  and $\mathbf{p}_{\beta}(k,n x)$ is as defined in \eqref{18.1.eq5}.\\
\begin{lemma}
Let $0<\alpha<\infty$, $| \beta|<1$ and $r\in \mathbb{N}$,
$$ S(r,\alpha,\beta) = \sum_{k=0 }^{\infty} \frac{1}{k!} (\alpha +
\beta k )^{k+r}e^{-(\alpha + \beta k)} $$ and
$$ (1-\beta)S(0,\alpha,\beta) = 1. $$
Then
$$ S(r,\alpha,\beta) = \sum_{k=0}^{\infty} \beta^k (\alpha + k\beta) S( r-1, \alpha+k\beta,\beta).$$
\end{lemma}
\textbf{Proof: } Notice that
\begin{eqnarray}\label{18.2.eq2}
 S(r,\alpha,\beta) &=&  \alpha \sum_{k=0 }^{\infty} \frac{1}{k!}  (\alpha + \beta k )^{k+r-1}e^{-(\alpha + \beta k)}
            +  \beta\sum_{k=0 }^{\infty} \frac{1}{k!} k (\alpha + \beta k )^{k+r-1}e^{-(\alpha + \beta k)}\nonumber\\
 &=& \alpha S(r-1,\alpha,\beta) + \beta\sum_{k=0 }^{\infty} \frac{1}{k!}     (\alpha + \beta +\beta k)^{k+r}e^{-(\alpha + \beta +\beta k)}\nonumber\\
 &=& \alpha S(r-1,\alpha,\beta) + \beta S(r,\alpha+\beta,\beta).
\end{eqnarray}
By a repeated use of \eqref{18.2.eq2}, the proof of the lemma is
archived.\\
Now when $|\beta | < 1$ , we have
\begin{eqnarray}
S(1,\alpha,\beta)&=&\sum_{k=0}^{\infty} \frac{\beta^k (\alpha +
k\beta)}{1-\beta}
=\frac{\alpha}{(1-\beta)^2}+\frac{\beta^2}{(1-\beta)^3};\\
S(2,\alpha,\beta) &=&  \sum_{k=0}^{\infty} \beta^k(\alpha+k
\beta)\left(\frac{\alpha+k\beta}{(1-\beta)^2}+\frac{\beta^2}{(1-\beta)^3}\right)
= \frac{\alpha^2}{(1-\beta)^3} +
\frac{3\alpha\beta^2}{(1-\beta)^4}+
\frac{\beta^3(1+2\beta)}{(1-\beta)^5}; \\
S(3,\alpha,\beta) &=& \sum_{k=0}^{\infty}\beta^k(\alpha+k \beta)
S(2,\alpha+k\beta,\beta) =\frac{\alpha ^3}{(1-\beta )^4}+\frac{6
\alpha ^2 \beta ^2}{(1-\beta )^5}+\frac{\alpha \beta^3 \left(4 +11
\beta\right)}{(1-\beta )^6}+\frac{\beta ^4+8 \beta ^5+6 \beta
^6}{(1-\beta )^7};\\
S(4,\alpha,\beta) &=& \sum_{k=0}^{\infty}\beta^k(\alpha+k \beta) S(3,\alpha+k\beta,\beta)\nonumber\\
&=& \frac{\alpha ^4}{(1-\beta )^5}+\frac{10 \alpha ^3 \beta
^2}{(1-\beta )^6}+\frac{5 \alpha ^2 \left(2 \beta ^3+7 \beta
^4\right)}{(1-\beta )^7}\nonumber\\
&&+\frac{5 \alpha  \left(\beta ^4+10 \beta ^5+10 \beta
^6\right)}{(1-\beta )^8}+\frac{\beta ^5+22 \beta ^6+58 \beta ^7+24
\beta ^8}{(1-\beta )^{9}}.
\end{eqnarray}
 \section{Estimation of Moments}
   We should note that, the first moment of the operators \eqref{18.1.eq3} gives  $J_n^{[\beta]}(t,x)=ax$, for some real constant $a$, but for the operators \eqref{18.2.eq1}, we have $P_n^{[\beta]}(t,x)=cx+b$ for some real constants $b$ and $c$. Due to this operators \eqref{18.1.eq3} and \eqref{18.2.eq1} have different approximation properties. We required following results to prove main results.
\begin{lemma}\label{18.2.lemma1}
The operators $P_{n}^{[\beta]}$ $n \geq 1$, defined by
\eqref{18.2.eq1} satisfy the following relations
 \begin{enumerate}
 \item $P_{n}^{[\beta]}(1,x)=1;$
 \item $\displaystyle P_{n}^{[\beta]}(t,x) = \frac{x}{(1-\beta)} + \frac{\beta}{n(1-\beta)^2};$
 \item $\displaystyle P_{n}^{[\beta]}(t^2,x) =  \frac{x^2}{(1-\beta )^2}+\frac{x (1+2 \beta )}{n (1-\beta )^3}+\frac{\beta  (1+2 \beta )}{n^2 (1-\beta )^4};$
 \item $\displaystyle P_{n}^{[\beta]}(t^3,x) = \frac{x^3}{(1-\beta )^3}+\frac{3 x^2 (1+\beta )}{n (1-\beta
)^4} +\frac{x \left(1+8 \beta +6 \beta ^2\right)}{n^2 (1-\beta
)^5}+\frac{\beta  \left(1+8 \beta +6 \beta ^2\right)}{n^3 (1-\beta
)^6};$
 \item $\displaystyle P_{n}^{[\beta]}(t^4,x) = \frac{x^4}{(1-\beta )^4}+\frac{2 x^3 (3+2 \beta )}{n (1-\beta )^5}+\frac{x^2 \left(7+26 \beta +12 \beta ^2\right)}{n^2 (1-\beta
 )^6}\\
~~~~~~~~~~~~~~~~~~~~~~~~~~~~ +\frac{x \left(1+22 \beta +58 \beta
^2+24 \beta ^3\right)}{n^3 (1-\beta )^7}+\frac{\beta  \left(1+22
\beta +58 \beta ^2+24 \beta ^3\right)}{n^4 (1-\beta )^8}$.
      \end{enumerate}
\end{lemma}
%\begin{proof}
 \textbf{Proof:}
By the relation (\ref{18.1.eq4}), it clear that $P_{n}^{[\beta]}(1,x)=1$.\\
 By the simple computation, we get
 \begin{eqnarray*}
 P_{n}^{[\beta]}(t,x) &=& (1-\beta)\sum_{k=0}^{\infty} ( n x + k\beta)^{k} \frac{e^{-(nx+k\beta)}}{(k)!}\frac{k}{n} \\
&=& \frac{(1-\beta)}{n} S(1,n x+\beta,\beta)= \frac{x}{(1-\beta)} + \frac{\beta}{n(1-\beta)}+\frac{\beta^2}{n(1-\beta)^2}
            = \frac{x}{(1-\beta)} + \frac{\beta}{n(1-\beta)^2};\\
P_{n}^{[\beta]}(t^2,x) &=& \frac{(1-\beta)}{n^2}\left\{ S(2,n x+2\beta,\beta) + S(1,n x+\beta,\beta)\right\}\\
&=& \frac{(1-\beta)}{n^2}\left\{\frac{(n x +
2\beta)^2}{(1-\beta)^3} + \frac{3(n x +
2\beta)\beta^2}{(1-\beta)^4}
                    + \frac{\beta^3(1+2\beta)}{(1-\beta)^5}+ \frac{(n x + \beta)}{(1-\beta)^2}+\frac{\beta^2}{(1-\beta)^3}\right\}\\
&=& \frac{x^2}{(1-\beta )^2}+\frac{x (1+2 \beta )}{n (1-\beta
)^3}+\frac{\beta  (1+2 \beta )}{n^2 (1-\beta )^4};
\end{eqnarray*}
\begin{eqnarray*}
P_{n}^{[\beta]}(t^3,x) &=& \frac{(1-\beta)}{n^3}\left\{S(3, n x +3\beta,\beta) + 3S(2,n x+2\beta,\beta) +S(1,n x+\beta,\beta) \right\}\\
&=&  \frac{(1-\beta)}{n^3}\left\{\frac{(n x + 3\beta) ^3}{(1-\beta )^4}+\frac{6 (n x + 3\beta) ^2 \beta ^2}{(1-\beta )^5}
+\frac{(n x + 3\beta)  \left(4 \beta ^3+11 \beta ^4\right)}{(1-\beta )^6}+\frac{\beta ^4+8 \beta ^5+6 \beta ^6}{(1-\beta )^7}\right. \\
&&\left.+ 3\left(\frac{(n x + 2\beta)^2}{(1-\beta)^3} + \frac{3(n x + 2\beta)\beta^2}{(1-\beta)^4}+ \frac{\beta^3(1+2\beta)}{(1-\beta)^5}\right)
+\frac{(n x + \beta)}{(1-\beta)^2}+\frac{\beta^2}{(1-\beta)^3} \right\}\\
&=& \frac{x^3}{(1-\beta )^3}+\frac{3 x^2 (1+\beta )}{n (1-\beta
)^4}
+\frac{x \left(1+8 \beta +6 \beta ^2\right)}{n^2 (1-\beta )^5}+\frac{\beta  \left(1+8 \beta +6 \beta ^2\right)}{n^3 (1-\beta )^6};\\
P_{n}^{[\beta]}(t^4,x)&=& \frac{(1-\beta)}{n^4}\left( S(4,n x
+4\beta,\beta)+ 6S(3,n x +3\beta,\beta)+7S(2,n x +2\beta,\beta)  +
S(1,n x +\beta,\beta)\right)\\
 &=& \frac{x^4}{(1-\beta )^4}+\frac{2 x^3 (3+2 \beta )}{n (1-\beta )^5}+\frac{x^2 \left(7+26 \beta +12 \beta ^2\right)}{n^2 (1-\beta
 )^6}\\
&& +\frac{x \left(1+22 \beta +58 \beta ^2+24 \beta ^3\right)}{n^3
(1-\beta )^7}+\frac{\beta  \left(1+22 \beta +58 \beta ^2+24 \beta
^3\right)}{n^4 (1-\beta )^8}.
\end{eqnarray*}
%\end{proof}
The proof of Lemma \ref{18.2.lemma1} is complete.
%\begin{lemma}\label{lemma.2.lem}
%Let $r \in \mathbb{N}$ be fixed number. Then there exist positive numerical coefficients $\lambda_{r,j,\beta}$, $1 \leq j \leq r$, depending only on $r$, $j$ and $\beta$ and $\varDelta_{\beta}$, depending only on $\beta$ such that
%$$ P_{n}^{[\beta]}(t^r,x)=\frac{1}{n^r(1-\beta)^r}\sum_{j=1}^r \frac{\lambda_{r,j,\beta}}{(1-\beta)^{j-1}}(nx)^j+\frac{\varDelta_{\beta}}{n^r(1-\beta)^{r+2}},$$
%for all $x \in [0,\infty)$ and $n \in \mathbb{N}$. Moreover, we have $\lambda_{1,1,\beta}=1$.
%\end{lemma}
%The proof follows by a mathematical induction argument.\\
	 We also introduce the $s$-th order central moment of the operator
$P_{n}^{[\beta]}$ , that is $P_{n}^{[\beta]}(\varphi_x^s,x)$,
where $\varphi_x(t) = t - x$, $(x, t) \in \mathbf{R}^+
\times\mathbf{R}^+$. On the basis of above lemma and by linearity
of operators (\ref{18.2.eq1}), by a straightforward calculation,
we obtain
\begin{lemma}\label{18.2.lemma4}
Let the operator $P_{n}^{[\beta]}$  be defined by relation as
(\ref{18.2.eq1}) and let $\varphi_x=t-x$ be given by
\begin{enumerate}
\item $\displaystyle P_{n}^{[\beta]}(\varphi_x,x) = \frac{x \beta}{(1- \beta)}+ \frac{\beta}{n(1-\beta)^2}$;
\item $\displaystyle  P_{n}^{[\beta]}(\varphi_x^2,x) = \frac{x^2 \beta ^2}{(1-\beta )^2}
+\frac{x \left(1+2 \beta ^2\right)}{n (1-\beta )^3}+\frac{\beta
(1+2 \beta )}{n^2 (1-\beta )^4};$
\item $\displaystyle P_{n}^{[\beta]}(\varphi_x^3,x) =\frac{x^3 \beta ^3}{(1-\beta )^3}+\frac{3 x^2 \beta  \left(1+\beta ^2\right)}{n (1-\beta )^4}
+\frac{x \left(1+5 \beta +3 \beta ^2+6 \beta ^3\right)}{n^2
(1-\beta )^5}+\frac{\beta  \left(1+8 \beta +6 \beta ^2\right)}{n^3
(1-\beta )^6};$
\item $\displaystyle P_{n}^{[\beta]}(\varphi_x^4,x) =\frac{x^4 \beta ^4}{(1-\beta )^4}+\frac{2 x^3 \beta ^2 \left(3+2 \beta ^2\right)}{n (1-\beta )^5}
 +\frac{x^2 \left(3+4 \beta +20 \beta
^2+6 \beta ^3+12 \beta ^4\right)}{n^2 (1-\beta )^6} \\
~~~~~~~~~~~~~~~~~~~~~~~~~~~~~+\frac{x \left(1+18 \beta +30 \beta
^2+32 \beta ^3+24 \beta ^4\right)}{n^3 (1-\beta )^7}+\frac{\beta
\left(1+22 \beta +58 \beta ^2+24 \beta ^3\right)}{n^4 (1-\beta
)^8} $.
\end{enumerate}
\end{lemma}
\begin{lemma}\label{18.2.lemma3}
Let $n > 1$ be a given number. For every $0<\beta<1$, one has
\begin{eqnarray*}
P_{n}^{[\beta]}(\varphi_x^2,x)\leq \frac{3+n}{n (1-\beta
)^4}\left(\phi(x)+\frac{1}{n}\right),
\end{eqnarray*}
where $\phi^2(x)=x(1+x)$, $x\in [0,\infty)$.
\end{lemma}
\textbf{Proof: }Since, $\max\{x,x^2\}\leq x+x^2$, $0<\beta<1$ and
$(1-\beta)^{-3} \leq  (1-\beta)^{-4}$, we have
\begin{eqnarray*}
P_{n}^{[\beta]}(\varphi_x^2,x)&=& (x^2 +x)\left(\frac{ \beta
^2}{(1-\beta )^2} +\frac{ \left(1+2 \beta ^2\right)}{n (1-\beta
)^3}\right)+\frac{\beta (1+2 \beta )}{n^2 (1-\beta )^4};\\
&\leq& \phi^2(x)\left(\frac{1+(2+n) \beta ^2-n \beta ^3}{n
(1-\beta )^3}\right) +\frac{\beta (1+2 \beta )}{n^2 (1-\beta )^4};\\
&\leq& \phi^2(x)\left(\frac{1+(2+n) \beta ^2}{n (1-\beta
)^3}\right) +\frac{\beta (1+2 \beta )}{n^2
(1-\beta )^4};\\
&\leq& \frac{1+(2+n) \beta ^2}{n (1-\beta
)^4}\left(\phi^2(x)+\frac{1}{n}\right),
\end{eqnarray*}
which is required.
\section{Approximation Properties}
 The convergence property of the
operators (\ref{18.2.eq1}) is proved in the following theorem:
\begin{theorem}\label{18.3.thm1}
Let $f$ be a continuous function on $[0,\infty)$ and $\beta_n \to
0$ as $n \to \infty$, then the sequence $P_{n}^{[\beta_n]}$
converges uniformly to f on $[a, b]$, where $0 \leq  a < b <
\infty$.
\end{theorem}
\textbf{Proof: }Since $P_{n}^{[\beta_n]}$ is a positive linear
operator for $0\leq  \beta_n<1$, it is
sufficient, by Korovkin's result \cite{korovkin1953convergence}, to verify the uniform convergence for test functions $f(t) = 1, t$ and $t^2$.\\
It is clear that $\displaystyle P_{n}^{[\beta_n]}(1,x) =1.$\\
Going to $f(t) = t $,
$$\lim_{n\to \infty}P_{n}^{[\beta_n]}(t,x)
         = \lim_{n\to \infty} \left[ \frac{x}{(1-\beta_n)} + \frac{\beta_n}{n(1-\beta_n)^2}\right]= x, \textrm { as } \beta_n \to 0.$$
Proceeding to the function $f(t) = t^2$, it can easily be shown
that
$$\lim_{n\to \infty}P_{n}^{[\beta_n]}(t^2,x) = \lim_{n\to \infty}\left[\frac{x^2}{(1-\beta_n )^2}
+\frac{x (1+2 \beta_n )}{n (1-\beta_n )^3}+\frac{\beta_n  (1+2
\beta_n )}{n^2 (1-\beta_n )^4}\right] = x^2, \textrm { as }
\beta_n \to 0.$$ The proof of theorem \ref{18.3.thm1} is complete.
\subsection{Local Approximation}
Let $C_B[0,\infty)$ be denote the set of all bounded
continuous real-valued functions on $[0,\infty)$. The space is
endowed with sup-norm $\|\cdot\|$, where $\displaystyle \|f\|=\sup_{x\in
[0,\infty)} |f(x)| $, $f\in C_B[0,\infty)$. In connection with the
estimation of the degree of approximation, the so
called moduli of smoothness play important role.\\
Further, let us consider the following $K$-functional:
$$K_2(f,\delta)=\inf_{g\in W^2} \{ \|f-g\|+ \delta\|g''\|\},$$
where $\delta > 0$ and $W^2 =\{g \in C_B[0,\infty): g'; g''\in
[0,\infty)\}$. By \cite[14, p. 177,Theorem
2.4]{devore1993constructive} there exists an absolute constant $C
> 0$ such that
\begin{equation}\label{18.2.eq3}
K_2(f,\delta) \leq C \omega_2(f,\sqrt{\delta}),
\end{equation}
where
$$\omega_2(f, \sqrt{\delta})=\sup_{0< h \leq \sqrt{\delta}} \sup_{x\in [0,\infty)}|f(x+2h) - 2f(x+h) +f(x)|$$
is the second order modulus of smoothness of $f\in
C_B[0,\infty)$.By
$$\omega(f,\delta)=\sup_{0<h\leq
\delta}\sup_{x\in[0,\infty)}|f(x+h)-f(x)|$$ we denote the usual
modulus of continuity of $f\in  C_B[0,\infty)$. In what follows we
shall use the notations $\phi(x) = \sqrt{x(1+x)}$ and $\delta^2(x)
=\phi^2(x)+\frac{1}{n}$, where $x \in [0,\infty)$ and $n\geq 1$.\\
Now, we establish local approximation theorems in connection with
the operators $P_n^{[\beta]}$ .
\begin{theorem}
Let $P_n^{[\beta]}$ , $n\in \mathbb{N}$, be given by
\eqref{18.2.eq1}. For every $f\in C_B[0,\infty)$, one has
\begin{eqnarray}
|{P}_n^{[\beta]}(f,x)-f(x)|\leq
C\omega_2\left(f,\frac{1}{2}\sqrt{\frac{1+(2+n) \beta ^2}{n
(1-\beta )^4}\delta^2(x)+\left(\frac{\beta (1+n x (1-\beta ))}{n
(1-\beta )^2}\right)^2}\right)+ \omega\left(f,\frac{\beta (1+n x
(1-\beta ))}{n (1-\beta )^2}\right).
\end{eqnarray}
\end{theorem}
\textbf{Proof: } Let us introduce the auxiliary operators
$\widetilde{P}_n^{[\beta]}$ defined by
\begin{equation}\label{18.2.eq4}
\widetilde{P}_n^{[\beta]}(f,x) = {P}_n^{[\beta]}(f,x) -
f\left(\frac{\beta +nx (1- \beta )}{n (1-\beta )^2}\right)+f(x),
\end{equation}
for $x\in [0,\infty)$. The operators $\widetilde{P}_n^{[\beta]}$
are linear. By Lemma \ref{18.2.lemma4}, we have
\begin{equation}\label{18.2.eq5}
\widetilde{P}_n^{[\beta]}(t-x, x) =0.
\end{equation}
Let $g\in W^2$. From Taylor's expansion
\begin{equation}
g(t) = g(x) + g'(x) (t-x) + \int_x^t (t-u) g''(u)du, ~~~t\in
[0,\infty).
\end{equation}
Applying the linear operator $\widetilde{P}_n^{[\beta]}$ and
taking in view \eqref{18.2.eq4} and \eqref{18.2.eq5}, we can write
\begin{eqnarray}\label{18.2.eq8}
|\widetilde{P}_n^{[\beta]}(g,x) -g(x)|&=&
|\widetilde{P}_n^{[\beta]}(g-g(x),x)|= \bigg|g'(x)
\widetilde{P}_n^{[\beta]}(t-x,x)+
\widetilde{P}_n^{[\beta]}\left(\int_x^{t}(t-u)g''(u)du,x\right)
\bigg|\notag\\
&\leq&\bigg|{P}_n^{[\beta]}\left(\int_x^{t}(t-u)g''(u)du,x\right)-\int_x^{\frac{\beta
+n x (1- \beta )}{n (1-\beta )^2}}\left(\frac{\beta +n x (1- \beta
)}{n (1-\beta )^2}-u\right)g''(u)du \bigg|\notag\\
&\leq&{P}_n^{[\beta]}\left(\bigg|\int_x^{t}(t-u)g''(u)du\bigg|,x\right)+\bigg|\int_x^{\frac{\beta
+n x (1- \beta )}{n (1-\beta )^2}}\left(\frac{\beta +n x (1- \beta
)}{n (1-\beta )^2}-u\right)|g''(u)|du \bigg|\notag\\
&\leq&{P}_n^{[\beta]}\left((t-x)^2\frac{\|g''\|}{2},x\right)+\|g''\|\bigg|\int_x^{\frac{\beta
+n x (1- \beta )}{n (1-\beta )^2}}\left(\frac{\beta +n x (1- \beta
)}{n (1-\beta )^2}-u\right)du \bigg|\notag\\
&\leq&\frac{\|g''\|}{2}{P}_n^{[\beta]}\left(\varphi_x^2,x\right)+\frac{\|g''\|}{2}\left(\frac{\beta
(1+n x (1-\beta ))}{n (1-\beta )^2}\right)^2\notag\\
&\leq&\|g''\|\left( \frac{1+(2+n) \beta ^2}{n (1-\beta
)^4}\delta^2(x)+\left(\frac{\beta (1+n x (1-\beta ))}{n (1-\beta
)^2}\right)^2\right).
\end{eqnarray}
 Let $f\in
C_B[0,\infty)$, Further on, taking in view that
\begin{equation}\label{18.2.eq6}|{P}_n^{[\beta]}(f-g,x)|\leq \|f-g\|,
~~~~~|\widetilde{P}_n^{[\beta]}(f-g,x)|\leq 3\|f-g\|\end{equation}
and by definition of modulus of continuity, we have
\begin{equation}\label{18.2.eq7}\bigg|f\left(\frac{\beta +nx (1- \beta )}{n (1-\beta
)^2}\right)-f(x)\bigg|\leq \omega\left(f,\frac{\beta  (1+n x
(1-\beta ))}{n (1-\beta )^2}\right).\end{equation} Now,
\eqref{18.2.eq8}, \eqref{18.2.eq6} and \eqref{18.2.eq7} imply
\begin{eqnarray*}
|{P}_n^{[\beta]}(f,x)-f(x)|&\leq&
|\widetilde{P}_n^{[\beta]}(f-g,x)-(f-g)(x)| + |
\widetilde{P}_n^{[\beta]}(g,x) -g(x)| + \bigg|f\left(\frac{\beta
+n x (1- \beta )}{n (1-\beta )^2}\right)-f(x)\bigg|\\
&\leq& 4\left(\|f-g\| +\frac{\|g''\|}{4}\left( \frac{1+(2+n) \beta
^2}{n (1-\beta )^4}\delta^2(x)+\left(\frac{\beta (1+n x (1-\beta
))}{n
(1-\beta )^2}\right)^2\right)\right)\\
&&+\omega\left(f,\frac{\beta (1+n x (1-\beta ))}{n (1-\beta
)^2}\right).
\end{eqnarray*}
Hence, taking infimum on the right hand side over all $g\in W^2$,
we get
\begin{eqnarray*}
|{P}_n^{[\beta]}(f,x)-f(x)|&\leq& 4 K_2\left(f,\frac{1}{4}\left(
\frac{1+(2+n) \beta ^2}{n (1-\beta
)^4}\delta^2(x)+\left(\frac{\beta (1+n x (1-\beta ))}{n (1-\beta
)^2}\right)^2\right)\right)+ \omega\left(f,\frac{\beta (1+n x
(1-\beta ))}{n (1-\beta )^2}\right).
\end{eqnarray*}
In view of (\ref{18.2.eq3}), we get
\begin{eqnarray*}
|{P}_n^{[\beta]}(f,x)-f(x)|&\leq&
C\omega_2\left(f,\frac{1}{2}\sqrt{\frac{1+(2+n) \beta ^2}{n
(1-\beta )^4}\delta^2(x)+\left(\frac{\beta (1+n x (1-\beta ))}{n
(1-\beta )^2}\right)^2}\right)+ \omega\left(f,\frac{\beta (1+n x
(1-\beta ))}{n (1-\beta )^2}\right).
\end{eqnarray*}
This completes the proof of the theorem.\\
We recall that a continuous function $f$ defined on $J$ is locally
$Lip$ $\alpha$ on $E$ ($0< \alpha \leq 1)$, if it satisfies the
condition
 \begin{equation}\label{18.3.eq1}
  |f(x)-f(y)| \leq M_f|x-y|^{\alpha},~~~(x,y) \in J \times E,
  \end{equation}
 where $M_f$ is a constant depending only on $f$.
\begin{theorem}
Let $P_{n}^{[\beta]}$, $n\in \mathbb{N}$, be given by
\eqref{18.2.eq1}, $0< \alpha\leq 1$  and $E$ be any subset of
$[0,\infty)$. If $f$ is locally $Lip$ $\alpha$  on $E$, then we
have
$$|P_{n}^{[\beta]}(f,x) -f(x) | \leq M_f \left(\left(\frac{1+(2+n) \beta ^2}{n (1-\beta
)^4}\delta^2(x)\right)^{\alpha/2}+ 2 d^{\alpha}(x,E) \right),
~~~x\geq 0,$$ where $d(x,E)$ is the distance between $x$ and $E$
defined as
$$d(x,E) = \inf \{ |x-y|: y\in E\}.$$
\end{theorem}
\textbf{Proof: }By using the continuity of $f$, it is obvious that
\eqref{18.3.eq1} holds for  any $x\geq 0$  and $y \in \bar{E}$ ,
$\bar{E}$ being  the closure in $\mathbb{R}$  of the set $E$. Let
$(x,x_0) \in [0,\infty) \times \bar{E}$  be such that  $|x-x_0| =
d(x,E)$.\\
On the other hand, we can write $|f(t) -f(x_0)| \leq |f(t)
-f(x_0)|+ |f(x_0) -f(x)|$  and applying the linear positive
operators $P_{n}^{[\beta_n]}$, we have
\begin{eqnarray*}
|P_{n}^{[\beta]}(f,x) - f(x)| &\leq&
|P_{n}^{[\beta]}(|f(t)-f(x_0)|,x) - f(x)| +|f(x_0) -f(x)|\\
&\leq & |P_{n}^{[\beta]}(M_f |t-x_0|^{\alpha} ,x)|
+M_f|x_0-x|^{\alpha}.
\end{eqnarray*}
Note that $P_{n}^{[\beta]}$  is positive, so it is monotone.\\
In the inequality $(A+B)^{\alpha}\leq A^{\alpha} + B^{\alpha}$ $(
A\geq 0, B\geq 0, 0< \alpha \leq 1)$, we put $A=|t-x|,$
$B=|x-x_0|$ and using Holder's inequality, we get
\begin{eqnarray*}
|P_{n}^{[\beta]}(f,x) - f(x)| &\leq&  M_f P_{n}^{[\beta]}(|(t-x)+(x-x_0)|^{\alpha},x) + M_f |x_0- x|^{\alpha}\\
&\leq& M_f \left(P_{n}^{[\beta]}(|t-x|^{\alpha},x) + |x-x_0|^{\alpha}\right) + M_f|x_0- x|^{\alpha}\\
&\leq& M_f \left(\left(P_{n}^{[\beta]}((t-x)^{2},x)\right)^{\alpha/2} +2 |x-x_0|^{\alpha}\right) \\
&\leq& M_f \left(\left(\frac{1+(2+n) \beta ^2}{n (1-\beta
)^4}\left(\phi(x)+\frac{1}{n}\right)\right)^{\alpha/2} +
2|x-x_0|^{\alpha}\right),
\end{eqnarray*}
which is required results.
\subsection{Rate of convergence} Let
$B_{x^2} [0;\infty)$ be the set of all functions $f$ defined on
$[0,\infty)$ satisfying the condition $\displaystyle |f(x)|\leq
M_f(1+x^2)$, where $M_f$ is a constant depending only on $f$. By
$C_{x^2}[0,\infty)$, we denote the subspace of all continuous
functions belonging to $B_{x^2}[0,\infty)$. Also, let
$C^*_{x^2}[0,\infty)$ be the subspace of all functions $f \in
C_{x^2}[0,\infty)$, for which $\displaystyle\lim_{x\to \infty}
\frac{f(x)}{1+x^2}$ is finite. The norm on $C_{x^2} [0\infty)$ is
$\displaystyle \|f\|_{x^2}=\sup_{x\in
[0,\infty)}\frac{|f(x)|}{1+x^2}$. For any positive $a$, by
$$\omega_a(f,\delta) = \sup_{|t-x|\leq \delta} \sup_{x,t\in [0,a]}|f(t)-f(x)|,$$
we denote the usual modulus of continuity of $f$ on the closed
interval $[0,a]$. We know that for a function $f\in
C_{x^2}[0,\infty)$, the modulus of continuity $\omega_a
(f,\delta)$ tends to zero. Now, we give a rate of convergence
theorem for the operator $P_n^{[\beta]}$:
\begin{theorem}
Let $f\in C_{x^2}[0,\infty)$ and let the operator
$P_{n}^{[\beta]}$ be defined as in (\ref{18.2.eq1}), where $0\leq
\beta < 1$ and $\omega_a(f,\delta)$ be its modulus of continuity
on the finite interval $[0,a]\subset [0,\infty)$, where $a > 0$.
Then for every $n\geq 1$,
\begin{eqnarray*}
\|P_{n}^{[\beta]}(f,\cdot) -f\| \leq \frac{(1+(2+n) \beta ^2)K}{n
(1-\beta )^4}+2\omega_{a+1}\left(f,\sqrt{\frac{(1+(2+n) \beta
^2)K}{n (1-\beta )^4}}\right),
\end{eqnarray*}
where $ K=6M_f(1+a^2)(1+a+a^2) $.
\end{theorem}
\textbf{Proof: } For $x\in [0,a]$ and $t>a+1$, since $t-x>1$, we
have
\begin{eqnarray}\label{18.2.4.eq1}
|f(t)-f(x)| &\leq& M_f ( 2+x^2+t^2)\notag\\
&\leq&M_f\left(2+3x^2 + 2(t-x)^2\right)\notag\\
&\leq& 6M_f (1+a^2) (t-x)^2.
\end{eqnarray}
For $x\in [0,a]$ and $t\leq a+1$, we have
\begin{equation}\label{18.2.4.eq2}
|f(t)-f(x)|\leq \omega_{a+1} (f, |t-x|) + \left(1+
\frac{|t-x|}{\delta}\right)\omega_{a+1}(f,\delta),
\end{equation}
with $\delta>0$. Form \eqref{18.2.4.eq1} and \eqref{18.2.4.eq2},
we have
\begin{equation}\label{18.2.4.eq3}
|f(t)-f(x)|\leq 6M_f(1+a^2)(t-x)^2+ \left(1+
\frac{|t-x|}{\delta}\right)\omega_{a+1}(f,\delta),
\end{equation}
for $x\in [0,a]$ and $t\geq 0$. Thus
\begin{equation}\label{18.2.4.eq4}
|P_n^{[\beta]}(f,x)-f(x)|\leq
6M_f(1+a^2)P_n^{[\beta]}((t-x)^2,x)+\omega_{a+1}(f,\delta)
\left(1+
\frac{1}{\delta}P_n^{[\beta]}\left((t-x)^2,x)\right)^{\frac{1}{2}}\right).
\end{equation}
Hence, by Schwarz's inequality and Lemma \ref{18.2.lemma3}, for
$0<\beta<1$ and $x\in [0,a]$
\begin{eqnarray}\label{18.2.4.eq5}
|P_n^{[\beta]}(f,x)-f(x)|&\leq& \frac{6M_f(1+a^2)(1+(2+n) \beta
^2)}{n (1-\beta
)^4}\left(\phi(x)+\frac{1}{n}\right)\notag\\
&&+\omega_{a+1}(f,\delta) \left(1+
\frac{1}{\delta}\sqrt{\frac{(1+(2+n) \beta ^2)}{n (1-\beta
)^4}\left(\phi(x)+\frac{1}{n}\right)}\right)\notag\\
&\leq& \frac{(1+(2+n) \beta ^2)K}{n (1-\beta
)^4}+\omega_{a+1}(f,\delta) \left(1+
\frac{1}{\delta}\sqrt{\frac{(1+(2+n) \beta ^2)K}{n (1-\beta
)^4}}\right).
\end{eqnarray}
By taking $\displaystyle \delta=\sqrt{\frac{(1+(2+n) \beta ^2)K}{n
(1-\beta )^4}}$, we get the assertion of our theorem.
\subsection{Weighted approximation}
Now, we shall discuss the weighted approximation theorem, where the
approximation formula holds true on the interval $[0,\infty)$.
\begin{theorem}
Let the operator $P_{n}^{[\beta_n]}$ be defined as in
(\ref{18.2.eq1}),  where $(\beta_n)_{n\geq 1},~ 0\leq \beta_n <
1$, satisfies $\displaystyle \lim_n\beta_n =0 $. For each $f\in
C_{x^2}^*[0,\infty)$, we have
$$\lim_{n\to \infty} \|P_n^{[\beta_n]}(f,\cdot)-f\|_{x^2}=0.$$
\end{theorem}
\textbf{Proof: } Using the theorem in \cite{Gadzhiev1976}, we see
that it is sufficient to verify the following three conditions
\begin{equation}\label{18.2.5.eq1}
\lim_{n\to \infty}\|P_n^{[\beta_n]}(t^v,\cdot)-x^v\|_{x^2}=0,
\textrm{ for } v=0,1,2,
\end{equation}
for every $x\in [0,\infty)$.\\
 Since $P_n^{[\beta_n]}(1,x)=1$, the first condition of (\ref{18.2.5.eq1}) is fulfilled for
 $v=0$.
 By Lemma \ref{18.2.lemma1} we have for $n \geq 1$
 \begin{eqnarray*}
\|P_n^{[\beta_n]}(t,x)-x\|_{x^2} &=&\sup_{x\in [0,\infty)}
\frac{|P_n^{[\beta_n]}(t,x)-x|}{1+x^2}\\
&\leq& \frac{\beta}{(1- \beta)}\sup_{x\in [0,\infty)}
\frac{x}{1+x^2}+ \frac{\beta}{n(1-\beta)^2}\leq \frac{\beta}{(1-
\beta)}+ \frac{\beta}{n(1-\beta)^2}
 \end{eqnarray*}
 and the second condition of (\ref{18.2.5.eq1}) holds for $v=1$ as $n \to \infty$ with $\beta_n\to
 0$.\\
 Similarly, we can write for $n\geq 1$
\begin{eqnarray*}
\|P_n^{[\beta_n]}(t^2,x)-x^2\|_{x^2} &\leq& \frac{(2-\beta_n )
\beta_n }{(1-\beta_n )^2}\sup_{x\in
[0,\infty)}\left(\frac{x^2}{1+x^2}\right)+\frac{ (1+2 \beta_n )}{n
(1-\beta_n )^3}\sup_{x\in [0,\infty)}\left(\frac{x^2}{1+x^2}
\right)+\frac{\beta_n  (1+2 \beta_n )}{n^2 (1-\beta_n )^4}\\
&\leq& \frac{(2-\beta_n ) \beta_n }{(1-\beta_n )^2}+\frac{ (1+2
\beta_n )}{n (1-\beta_n )^3}+\frac{\beta_n (1+2 \beta_n )}{n^2
(1-\beta_n )^4},
 \end{eqnarray*}
which implies that
$$\lim_{n\to \infty} \|P_n^{[\beta_n]}(t^2,x)-x^2\|_{x^2}=0 \textrm{ with } \beta_n \to \infty.$$
Thus the proof is completed.

\subsection{Asymptotic Formula} In order to present our asymptotic
formula, we need the following lemma.

\begin{lemma}\label{14.3.lemma1} Let $ P_{n}^{[\beta]}$ be defined as \eqref{18.2.eq1}. In addition, $0<\beta<1$, then
$$ P_{n}^{[\beta]}((t-x)^4,x)\leq \frac{267(x + x^2+x^3+x^4) }{n^4(1-\beta)^8}.$$
\end{lemma}
\textbf{Proof:} Since $\max\{x, x^2,x^3,x^4\} \leq x +
x^2+x^3+x^4$, $(1-\beta)^2 \leq 1$ and $(1-\beta)^{-i} \leq
(1-\beta)^{-(i+1)}$, for $i\in \mathbb{N}$,
\begin{eqnarray*}
P_{n}^{[\beta]}((t-x)^4,x) &\leq&\frac{x^4}{(1-\beta )^4}+\frac{
10x^3}{n (1-\beta )^5}
 +\frac{45x^2}{n^2 (1-\beta )^6} +\frac{105x}{n^3 (1-\beta )^7}+\frac{105}{n^4 (1-\beta
)^8}\\
&\leq& \frac{267(x + x^2+x^3+x^4) }{n^4(1-\beta)^8},
\end{eqnarray*}
we obtain our claim inequality.\\
Notice that, $\displaystyle \lim_{n\to \infty} n^2
P_{n}^{[\beta_n]}((t-x)^4,x) = 3x^2$ with $\beta_n\to 0$.
\begin{theorem} Let $f,f',f'' \in  C[0,\infty)$ and let the operator $P_{n}^{[\beta_n]}$ be defined as in (\ref{18.2.eq1}),
 where $(\beta_n)_{n\geq 1},~ 0\leq \beta_n < 1$, satisfies
 $\displaystyle \lim_n\beta_n =0 $. then
$$\lim_{n\to \infty} n\left( P_{n}^{[\beta_n]} (f,x) -f(x) \right) = \frac{x}{2}f''(x),~~ \forall ~~x>0. $$
\end{theorem}
\textbf{Proof: } Let $f,f',f'' \in  C[0,\infty)$ and $x \in
[0,\infty)$ be fixed. By the Taylor formula, we have
\begin{equation}\label{14.3.eq1}f(t) = f(x) + f'(x) (t-x) + \frac{1}{2}f''(x) (t-x)^2 + r(t;x)(t-x)^2,\end{equation}
where $r(\cdot; x)$ is the Peano form of the remainder and
$\displaystyle \lim_{t \to x} r(t;x) =0$.\\
We apply $P_{n}^{[\beta_n]}$ to equation (\ref{14.3.eq1}), we get
\begin{eqnarray*}P_{n}^{[\beta_n]}(f,x)- f(x)&=& f'(x) P_{n}^{[\beta_n]}((t-x),x)
            + \frac{1}{2}f''(x) P_{n}^{[\beta_n]}((t-x)^2,x) + P_{n}^{[\beta_n]}(r(t ;x)(t-x)^2,x)\\
&=&  f'(x)\left[\frac{x \beta_n}{(1- \beta_n)}+ \frac{\beta_n}{n(1-\beta_n)^2} \right]+ P_{n}^{[\beta_n]}(r(t ;x)(t-x)^2,x) \\
&&+  \frac{f''(x)}{2}\left[\frac{x^2 \beta_n ^2}{(1-\beta_n )^2}
+\frac{x \left(1+2 \beta_n ^2\right)}{n (1-\beta_n
)^3}+\frac{\beta_n (1+2 \beta_n )}{n^2 (1-\beta_n
)^4}\right].\end{eqnarray*} In the second term
$P_{n}^{[\beta_n]}(r(t ;x)(t-x)^2,x)$ applying the Cauchy-Schwartz
inequality, we have
\begin{equation}\label{14.3.eq2}
0\leq |n P_{n}^{[\beta_n]}(r(t ;x)(t-x)^2,x)| \leq
\sqrt{n^2P_{n}^{[\beta_n]}((t-x)^4,x)}\sqrt{P_{n}^{[\beta_n]}(r(t
;x)^2,x)}.\end{equation} 
Observe that $r^2(x,x)=0$ and $r^2(\cdot,x)\in C_{x^2}^*[0,\infty)$. Then, it follows that
\begin{equation}\label{14.3.eq3}\lim_{t\to x} nP_{n}^{[\beta_n]}(r(t,x)^2,x) =r^2(x,x)=0\end{equation} uniformly with respect to $x\in [0,A]$ for any $A>0$ .\\
On the basis of (\ref{14.3.eq2}), (\ref{14.3.eq3}) and Lemma
\ref{14.3.lemma1} , we get
\begin{eqnarray*}
\lim_{n\to \infty}n\left(P_{n}^{[\beta_n]}(f,x)- f(x)\right)&=&
f'(x)\lim_{n\to \infty}n\left[\frac{x \beta_n}{(1- \beta_n)}+
\frac{\beta_n}{n(1-\beta_n)^2}\right]\\
&&+ \frac{ f''(x)}{2}\lim_{n\to \infty}  n\left[\frac{x^2 \beta_n
^2}{(1-\beta_n )^2} +\frac{x \left(1+2 \beta_n ^2\right)}{n
(1-\beta_n )^3}+\frac{\beta_n (1+2 \beta_n )}{n^2 (1-\beta_n
)^4}\right]\\
&=& \frac{x}{2}f''(x) \textrm{ with } \beta_n \to 0,
\end{eqnarray*}
which completes the proof.

\section*{Acknowledgements}
The authors are thankful to the referees for valuable suggestions, leading to an overall improvement in the paper.


\begin{thebibliography}{10}
	\expandafter\ifx\csname url\endcsname\relax
	\def\url#1{\texttt{#1}}\fi
	\expandafter\ifx\csname urlprefix\endcsname\relax\def\urlprefix{URL }\fi
	\expandafter\ifx\csname href\endcsname\relax
	\def\href#1#2{#2} \def\path#1{#1}\fi
	
	\bibitem{consul1973generalization}
	P.~C. Consul, G.~C. Jain, A generalization of the poisson distribution,
	Technometrics 15~(4) (1973) 791--799.
	
	\bibitem{consul1973some}
	P.~C. Consul, G.~Jain, On some interesting properties of the generalized
	poisson distribution, Biometrische Zeitschrift 15~(7) (1973) 495--500.
	
	\bibitem{jain1975linear}
	G.~Jain, A linear function poisson distribution, Biometrische Zeitschrift
	17~(8) (1975) 501--506.
	
	\bibitem{Berstien1912Demo}
	S.~N. Berstien, D\'{e}monstration du th\'{e}or\'{e}me de \textsc{W}eierstrass
	fond\'{e}e sur le calcul de probabilit\'{e}s, Commun. Soc. Math. Kharkow
	13~(2) (1912-1913) 1--2.
	
	\bibitem{szasz1950generalization}
	O.~Sz{\'a}sz, Generalization of \textrm{S. B}ernsteins polynomials to the
	infinite interval, J. Res. Natl. Bur. Stand. 45 (1950) 239--245.
	
	\bibitem{meyer1960bernsteinsche}
	W.~Meyer-K{\"o}nig, K.~Zeller, Bernsteinsche potenzreihen, Studia Mathematica
	19~(1) (1960) 89--94.
	
	\bibitem{cheney1964generalization}
	E.~Cheney, A.~Sharma, On a generalization of \textrm{B}ernstein polynomials,
	riv, Mat. Univ. Parma (2) 5 (1964) 77--84.
	
	\bibitem{stancu1968approximation}
	D.~Stancu, Approximation of functions by a new class of linear polynomial
	operators, Rev. Roumaine Math. Pures Appl. 13~(8) (1968) 1173--1194.
	
	\bibitem{mirakyan1941approximation}
	G.~Mirakyan, Approximation des fonctions continues au moyen polyn{\^o}mes de la
	forme, Dokl. Akad. Nauk. SSSR 31 (1941) 201--205.
	
	\bibitem{rempulska2001approximation}
	L.~Rempulska, Z.~Walczak, Approximation properties of certain modified
	szasz-mirakyan operators, Le Matematiche 55~(1) (2001) 121--132.
	
	\bibitem{abel2007asymptotic}
	U.~Abel, M.~Ivan, X.-M. Zeng, Asymptotic expansion for szasz-mirakyan
	operators, in: Numerical Analysis and Applied Mathematics, Vol. 936, 2007,
	pp. 779--782.
	
	\bibitem{may1976saturation}
	C.~May, Saturation and inverse theorems for combinations of a class of
	exponential-type operators, Canad. J. Math 28~(6) (1976) 1224--1250.
	
	\bibitem{jain1972approximation}
	G.~C. Jain, Approximation of functions by a new class of linear operators,
	Journal of the Australian Mathematical Society 13~(3) (1972) 271--276.
	
	\bibitem{agratini2013approximation}
	O.~Agratini, Approximation properties of a class of linear operators,
	Mathematical Methods in the Applied Science 36~(17) (2013) 2353--2358.
	
	\bibitem{umar1985approximation}
	S.~Umar, Q.~Razi, Approximation of function by a generalized \textsc{S}zasz
	operators, Communications de la Facult{\'e} Des Sciences de L'Universit{\'e}
	D'Ankara: Math{\'e}matique 34 (1985) 45--52.
	
	\bibitem{tarabie2012jain}
	S.~Tarabie, On \textsc{J}ain-\textsc{B}eta linear operators, Applied
	Mathematics \& Information Sciences 6~(2) (2012) 213--216.
	
	\bibitem{mishrasome2013}
	V.~N. Mishra, P.~Patel, Some approximation properties of modified
	\textsc{J}ain-\textsc{B}eta operators, Journal of Calculus of Variations 2013
	(2013) 8 pages.
	
	\bibitem{Patelmishra20131}
	P.~Patel, V.~N. Mishra, Jain-\textrm{B}askakov operators and its different
	generalization, Acta Mathematica Vietnamica 40 (2015) 715--733.
	
	\bibitem{agratini2014approximation}
	O.~Agratini, On an approximation process of integral type, Applied Mathematics
	and Computation 236 (2014) 195--201.
	
	\bibitem{farcasapproximation}
	A.~FARCAS, An approximation property of the generalized \textrm{J}ains
	operators of two variables, Mathematical Sciences And Applications E-Notes
	1~(2) (2013) 158--164.
	
	\bibitem{lupas1995approximation}
	A.~Lupas, The approximation by some positive linear operators, in: Proceedings
	of the International Dortmund Meeting on Approximation Theory, Akademie
	Verlag, Berlin, 1995, pp. 201--229.
	
	\bibitem{patel2015onnew}
	P.~Patel, V.~N. Mishra, On new class of linear and positive operators,
	Bollettino dell'Unione Matematica Italiana 8~(2) (2015) 81--96.
	
	\bibitem{korovkin1953convergence}
	P.~Korovkin, On convergence of linear positive operators in the space of
	continuous functions, in: Dokl. Akad. Nauk SSSR, Vol.~90, 1953, pp. 961--964.
	
	\bibitem{devore1993constructive}
	R.~A. DeVore, G.~G. Lorentz, Constructive approximation, Vol. 303, Springer
	Verlag, 1993.
	
	\bibitem{Gadzhiev1976}
	A.~D. Gadzhiev, Theorems of the type of \textrm{P. P.} \textrm{K}orovkin type
	theorems, Math. Zametki 20 (5) (1976) 781--786, \text{Math Notes} 20 (5-6)
	(1976) 996-998 (English Translation).
	
\end{thebibliography}
\end{document}